\documentclass[11pt]{amsart}

\usepackage{amsthm}
\usepackage{amsfonts}
\usepackage[mathcal]{euscript}
\usepackage{amssymb}
\usepackage{amsrefs} 
\usepackage{graphicx}

\input xy 
\xyoption{all} 
\theoremstyle{plain}
\newtheorem{proposition}{Proposition}
\newtheorem{theorem}{Theorem}
\newtheorem{remark}{Remark}
\newtheorem{lemma}{Lemma}
\newtheorem{definition}{Definition}

\theoremstyle{definition}
\newtheorem{example}{Example}

\newcommand{\ra}{\rightarrow}

\newcommand{\bE}{\mathbb{E}}
\newcommand{\bB}{\mathbb{B}}

\title[Khovanov homology and diagonalisable Frobenius algebras]{Khovanov homology and diagonalisable Frobenius algebras}

\author{ Paul
  Turner}
\thanks{Supported by NCCR SwissMAP}
\address{
{\sc Paul Turner:} Section de math\'ematiques,  
Universit\'e de Gen\`eve, 2-4 rue du Li\`evre, CH-1211, Geneva}
\email{Paul.Turner@unige.ch}

\begin{document}

\vspace*{-1cm}

\begin{abstract}
We give a short elementary proof that a Khovanov-type link homology constructed from a diagonalisable Frobenius algebra is degenerate.
\end{abstract}

\maketitle

A central ingredient in  the construction of Khovanov homology is a particular rank two, commutative Frobenius algebra. Some other rank two Frobenius algebras, but not all, also give rise to Khovanov-type link homology theories, the best-known of these being the rational theory defined by Lee (\cite{0201105}) and the mod 2 theory defined by Bar-Natan (\cite{0410495}, see also \cite{0411225}). Both of these are ``degenerate'' in the sense that the homology has total rank depending only on the number of components and the degrees of the generators are a function of linking numbers among the components. One may also work over different rings (considering integral and mod $p$ Lee theory or Bar-Natan theory, for example) and the key to degeneracy lies in the fact that the Frobenius algebra is diagonalisable. Indeed the following result is considered to be ``well-known'':

\vspace{1.5mm}

\begin{theorem}\label{thm:main}
Let $A$ be a diagonalisable Frobenius algebra defined over a commutative ring $R$. Let $H^*(-)$ be a Khovanov-type link homology constructed from $A$. Then for any link $L$ with $|L|$ components, the homology $H^*(L)$ is a free module of rank $2^{|L|}$. Moreover there is a canonical basis corresponding to the set of orientations of $L$. 
\end{theorem}

\vspace{1.5mm}

However, the ``well-known proof'' of this result - which goes ``extend Lee's proof'' -  is flawed beyond rational Lee theory. At a crucial point a ``Hodge theory'' argument is used and this requires the use of a positive definite bilinear form. This argument does not, in fact, naively extend to the general case. I am grateful to Robert Lipshitz and Sucharit Sarkar who brought this to my attention and indirectly to Peter Kronheimer and Tom Mrowka who brought it to theirs. For Lee theory there other methods are available (see \cite{0606542}, for example) but it appears there is no proof that covers all cases of diagonalisable Frobenius algebras at once.

The main goal of this short note is to provide a simple argument that works in all cases. 

\vspace{4mm}

\section{The poset of enhanced states}
We will adopt an approach to Khovanov homology which is closely related to the one described by Viro in \cite{0202199}. Instead of organising the Kauffman states of a diagram into a cube we consider enhanced states which we will organise into a poset. 

The construction of the Khovanov-type homology of a link based on a Frobenius algebra $A$ is well known: from a link diagram assemble the set of Kauffman states of the diagram into a hyper-cube; associate to a Kauffman state $s$ (at a vertex of the hyper-cube) a free module $V_s$ being a tensor power of the Frobenius algebra; and finally, construct a cochain complex $C^*$ by taking direct sums of the $V_s$ and by using the Frobenius algebra multiplication and co-multiplication along the edges of the hyper-cube to define a differential - not forgetting that the edges need some signs to make things work. The homology of this complex is the homology of the link. 

It will be convenient for us to use the language of posets rather than ``cubes''. Recall that in a poset with partial order $\leq$ we say $y$ {\em covers} $x$, denoted $x\prec y$, if $x\leq y$ and for any $x\leq z\leq y$ we have either $z=x$ or $z=y$. Also recall that the {\em Hasse diagram} of a poset is the graph with vertices the elements of the poset with one edge for each covering relation. The  {\em Boolean lattice} on a set $S$, denoted $\bB(S)$, is the poset of subsets in $S$ ordered by inclusion and its Hasse diagram is a hyper-cube. Indeed, the Khovanov cube for a link diagram, is based on the Boolean lattice $\bB$ of subsets of the set of crossings. The labelling by modules may be expressed as a functor from the poset $\bB$, viewed as a category, to the category of $R$-modules, but in this paper we will simple continue to refer to ``labelling by modules''.

In an early paper on Khovanov homology (\cite{0202199}) Viro expresses the construction of the complex a little differently and we will paraphrase his approach here. Viro's approach was taken up in some subsequent work, notably that of Shumakovitch (\cite{0405474}). 

We suppose that the Frobenius algebra used in the construction has basis $\{a, b\}$ (at this stage not necessarily a diagonal basis) and instead of Kauffman states we will consider the set of enhanced states: an   {\em enhanced state} $(s,x)$ consists of a Kauffman state $s$ together with a labelling $x$ of its component circles by basis elements, that is to say a function $x\colon \{\text{circles in $s$}\} \ra \{a,b\}$. 

Given a Kauffman state $s$, the set of possible functions $x$ provides a basis for the (free) module $V_s$ and we may write
$$
V_s \cong \bigoplus_{x} R.
$$
Given an edge in the cube, which corresponds to a covering relation $s\prec s^\prime$, the original construction associates a map $V_s \ra V_{s^\prime}$, which, by using the above we may consider as a map $ \bigoplus_{x} R \ra \bigoplus_{x^\prime} R$. Such a map may be expressed in terms of its matrix components, which we denote $[x: x^\prime]\in R$. By definition, unless $s\prec s^\prime$ the numbers $[x: x^\prime]$ are zero. The complex which computes the homology may now be expressed as a direct sum over all enhanced states $(s,x)$ with differential determined by its matrix components $[x: x^\prime]$.

The enhanced states may be arranged into a poset in the following way.  

\vspace{1.5mm}

\begin{definition}
Given a link diagram, its {\em poset of enhanced states}, denoted $\bE$, is the poset whose elements are the enhanced states and for which the covering relation is defined by $(s,x) \prec (s^\prime,x^\prime)$ if and only if $[x: x^\prime]\neq 0$. \end{definition}

\vspace{1.5mm}

The partial order is the induced one: $(s,x) \leq (s^\prime,x^\prime)$ if either $(s,x) = (s^\prime,x^\prime)$ or there is a chain of covering relations from $(s,x)$ to $(s^\prime,x^\prime)$. There is a poset map $\bE \ra \bB$ defined by $(s,x) \mapsto s$ and it can be useful to visualise the Hasse diagram of the poset of enhanced states as ``living over'' the boolean lattice (in the sense of bundles).

Just as the Boolean lattice of a diagram is labelled with modules and maps we can consider the {\em labelled} poset of enhanced states, which associates to each element the ground ring $R$ and to a covering relation $(s,x) \prec (s^\prime,x^\prime)$ the map which is multiplication by $[x: x^\prime]$.

\vspace{1.5mm}

\begin{example}
Consider the unknot as presented by the diagram shown here: \raisebox{-2.5mm}{ \includegraphics[width=0.1\linewidth]{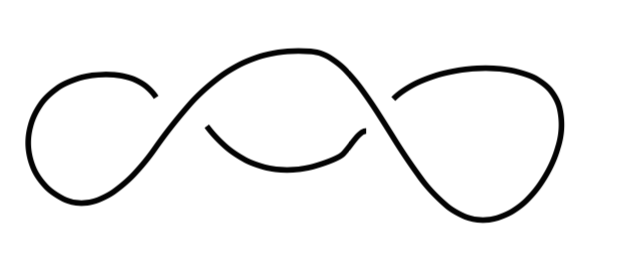}}.  We take the Frobenius algebra defining Lee theory which has basis ${a,b}$ with multiplication $a.a=a$, $a.b=b.a=0$, $b.b=b$ and co-multiplication  $\Delta(a) = 2a\otimes a$, $\Delta(b) = -2 b\otimes b$. (This is not quite the basis used by Lee: our $a$ is half hers and our $b$ is minus one half hers.) Here is a picture of the labelled poset of enhanced states lying above the Boolean lattice of Kauffman states. Any unlabelled edges in $\bE$ have matrix component 1. The edges labelled $-1$ get the minus sign from the sign assignment of the edge they lie above.

\vspace{-5pt}
\begin{center}
\includegraphics[width=0.8\linewidth]{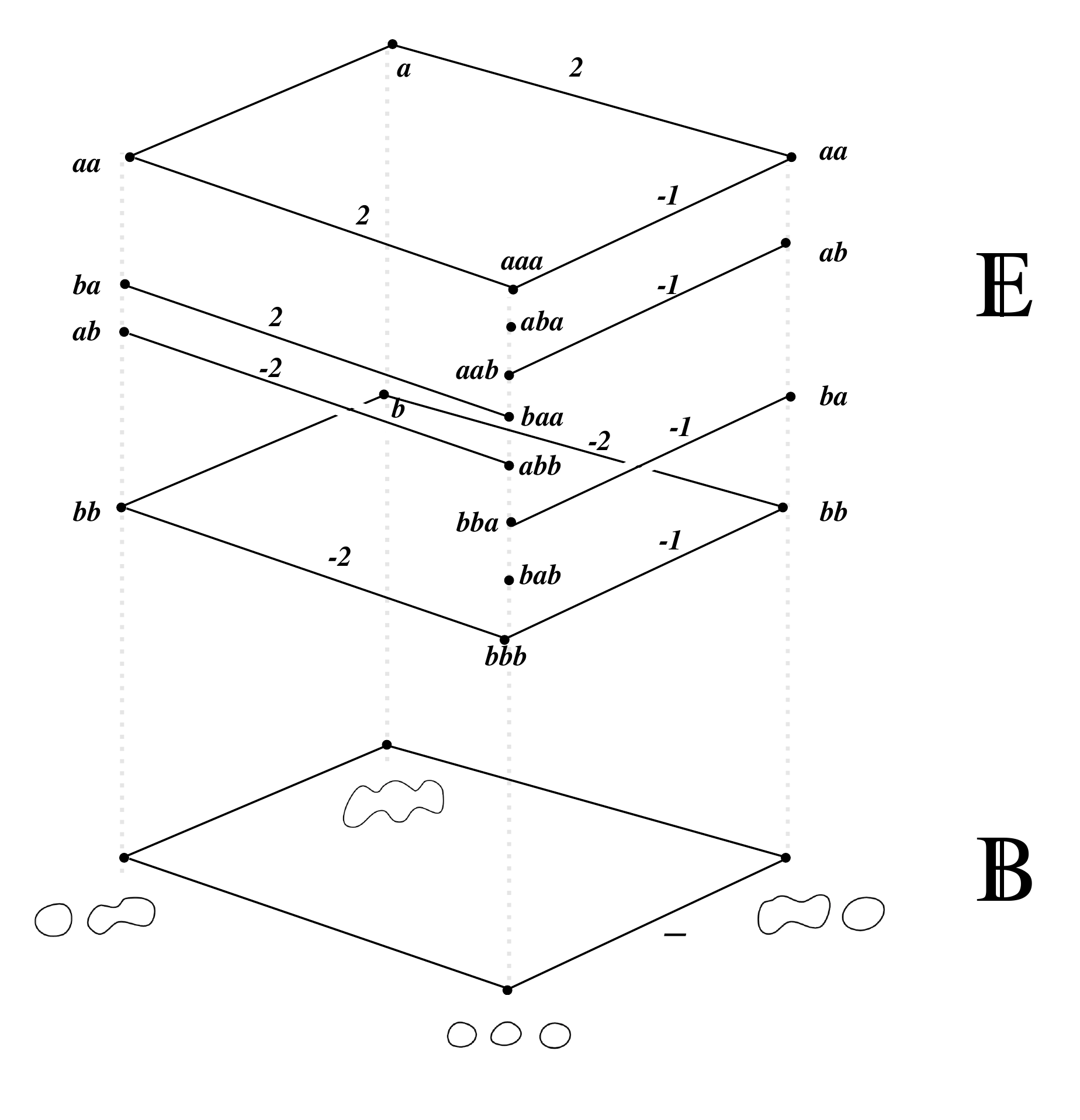}
\end{center}
\vspace{-15pt}
\end{example}

Notice that in this example the poset of enhanced states decomposes into a disjoint union of Boolean lattices: two of rank 2 (squares), four of rank 1 (lines) and two of rank 0 (points). As we will see in the next section this behaviour is typical of what happens for any link when the Frobenius algebra has multiplication of this form and  this is the key to proving the main result.

\vspace{4mm}

\section{Diagonalisable Frobenius Algebras and the poset of enhanced states}
The Frobenius algebra that defines Lee theory (over a ring in which 2 is invertible) and Bar-Natan theory (always) have a special property: they are diagonalisable in the following sense.

\vspace{1.5mm}

\begin{definition}
Let $R$ be a commutative ring with unit and $A$ a commutative Frobenius algebra over $R$ of rank two (as an $R$-module it is free of rank two). We say that $A$ is {\em diagonalisable} if there is a basis $\{a,b\}$ such that $a^2=a$, $b^2=b$ and $ab=0$. 
\end{definition}

\vspace{1.5mm}

In this section we show that if the Frobenius algebra used to define a link homology is diagonalisable, then the poset of enhanced states  decomposes into a disjoint union of Boolean lattices. 

\vspace{1.5mm}

\begin{lemma} \label{lem:coprod} Let $A$ be a diagonalisable Frobenius algebra with basis $\{a,b\}$. Writing  $\epsilon$ for the Frobenius algebra co-unit and $\Delta$ for its co-multiplication, the elements 
$\epsilon(a)$ and $\epsilon(b)$ are invertible in $R$ and 
$$
\Delta(a) = \epsilon(a)^{-1}a\otimes a \;\;\;\;\;\; \text{ and }  \;\;\;\;\;\;  \Delta(b) = \epsilon(b)^{-1}b\otimes b.
$$
\end{lemma}

\vspace{1.5mm}

\begin{proof}
Writing $\langle - , -\rangle$ for the Frobenius algebra inner product we have $\langle a, b\rangle = \epsilon(ab)= \epsilon(0)= 0$ and $\langle a, a\rangle = \epsilon(a)$ and thus by the non-degeneracy of $\langle -, -\rangle$ we see that $\epsilon(a)\neq 0$. Similarly we see $\epsilon(b)\neq 0$. 

Let us write $\Delta(a) = \alpha a\otimes a + \beta a\otimes b + \gamma b\otimes a + \delta b\otimes b$.  In any Frobenius algebra the map $ (1 \otimes \langle -, -\rangle) \circ (\Delta \otimes 1) \colon A\otimes A \ra A$ is precisely the multiplication map.  Feeding $a\otimes b$ into this gives
$$
\epsilon(b) \beta a + \epsilon(b) \delta b = 0.
$$
Since $\epsilon(b)\neq 0$ it follows that $\beta=\gamma=0$. Feeding $a\otimes a$ to the same equality gives
$$
\epsilon(a) \alpha a + \epsilon(a) \gamma b = a 
$$ 
which, since  $\epsilon(a)\neq 0$, gives $\gamma = 0$ and $\epsilon(a) \alpha = 1$. Thus, $\epsilon(a)$ is invertible (with inverse $\alpha$) and $\Delta(a) = \epsilon(a)^{-1}a\otimes a$.
A similar argument holds of $\epsilon(\beta)$ and $\Delta(b)$.

\end{proof}

From now on we assume we are working with a diagonalisable Frobenius algebra with basis $\{a,b\}$. 

In the process of smoothing crossings to obtain Kauffman states, the information about the position of crossings is lost. For convenience we wish to record the location of crossing information in the form of arcs as illustrated here:
 
\begin{center} 
\includegraphics[width=0.5\linewidth]{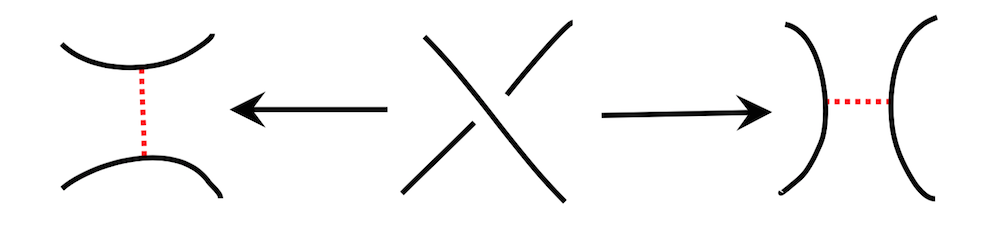}
\end{center}

Given an enhanced state such an arc may bridge sections of circle with the same label ($a$ or $b$) or with different labels.

\vspace{1.5mm}

\begin{definition}
Given an enhanced state $(s,x)$ we define its {\em arc-set} to be the set of arcs which connect sections of circle with the same label. The {\em arc-resolution} of an enhanced state $(s,x)$ consists of the union of its arc-set and the circles in the Kauffman state. 
\end{definition}

\vspace{1.5mm}

\begin{center} 
\includegraphics[width=0.7\linewidth]{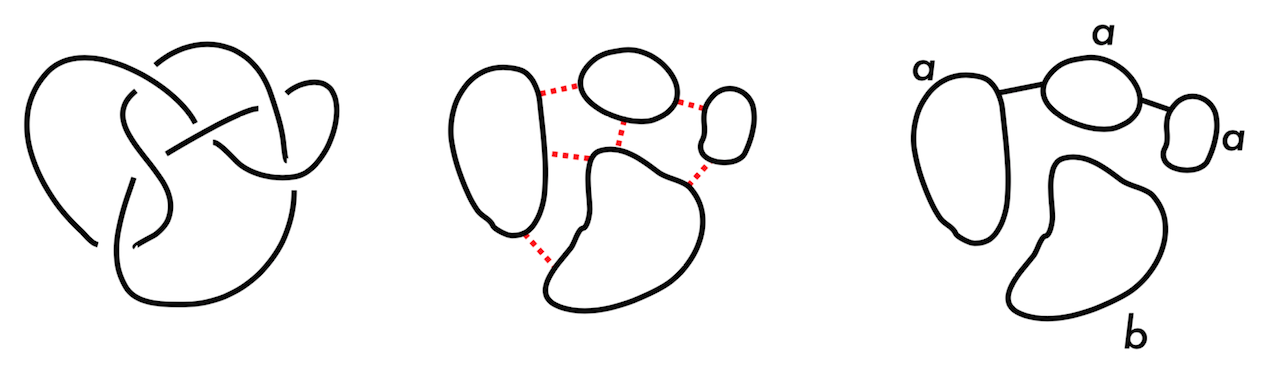}
\end{center}

Given two enhanced states 
$(s,x)$ and$ (s^\prime,x^\prime)$ 
we say that their arc-sets {\em correspond} if they are associated to the same set of crossings. In such a situation there is a natural bijection between the connected components of the associated arc-resolutions. If this bijection respects the labellings $x$ and $x^\prime$, in other words the following diagram commutes
$$
\xymatrix{
\{\text{connected components of $s$}\}   \ar[rr]^{\text{bijection}} \ar[dr]_{x} && \{\text{connected components of $s^\prime$}\}   \ar[dl]^{x^\prime}\\
& \{a,b\}
}
$$
then we say that the labellings {\em correspond}.

\vspace{1.5mm}

\begin{proposition}\label{prop:enhancedcompare}
Let $A$ be a diagonalisable Frobenius algebra and $\bE$ the poset of enhanced states of a link diagram based on $A$.  If $(s,x) \leq (s^\prime,x^\prime)$ in $ \bE$ then the arc-sets and labellings of $(s,x)$ and $(s^\prime,x^\prime)$ correspond.
\end{proposition}

\vspace{1.5mm}

\begin{proof}

It suffices to consider the case $(s,x)\prec (s^\prime,x^\prime)$ in $\bE$. We have $s\prec s^\prime$ in $\bB$ so at one crossing a 0-resolution is changed into a 1-resolution:

\begin{center} 
\includegraphics[width=0.4\linewidth]{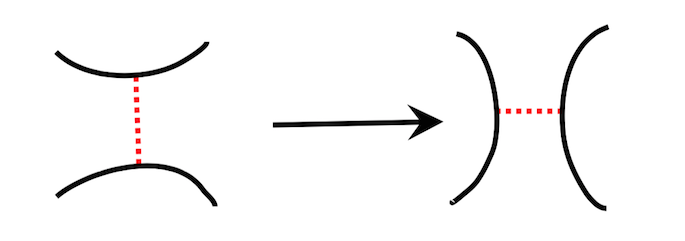}
\end{center}

A priori the dashed arcs shown may or may not be in the arc sets. There are two possibilities for connecting up the strands:

 \begin{center} 
(i)  \raisebox{-15mm}{\includegraphics[width=0.3\linewidth]{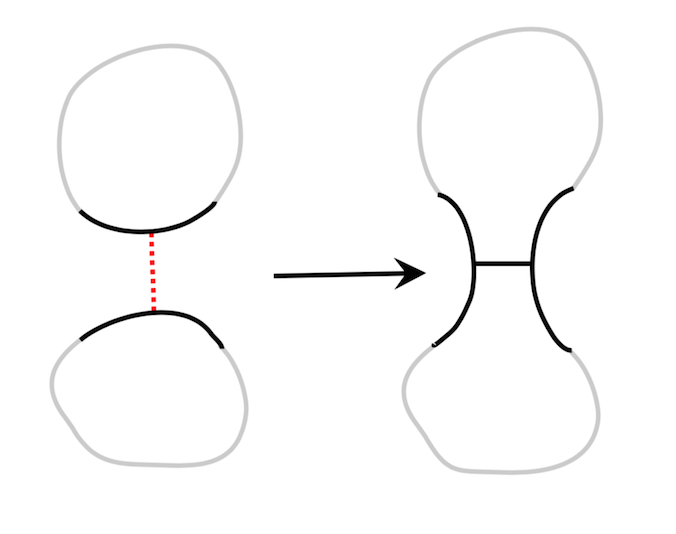}} \hspace*{8mm}
(ii) \raisebox{-8mm}{\includegraphics[width=0.5\linewidth]{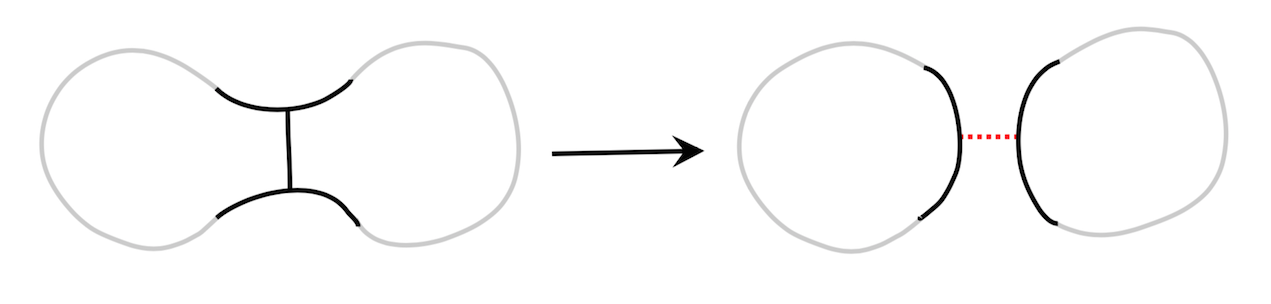}}
\end{center}

In case (i), the arc on the right is always present in the arc-resolution since it links parts of the same circle. Using the definition of multiplication and the fact that $[x:x^\prime]\neq 0$ the circles on the left must have identical labels (the same as the circle on the right) and so the (dashed) arc is, in fact, in the arc-resolution on the left. For case (ii) the arc on the left is similarly always present, and using the definition of co-multiplication and the fact that $[x:x^\prime]\neq 0$ we see that the circles on the right have identical labels  (the same as the circle on the left) and the (dashed) arc is once again in the arc-resolution on the right. Summarising, we have the following possibilities for enhanced states:

\begin{center} 
\includegraphics[width=0.6\linewidth]{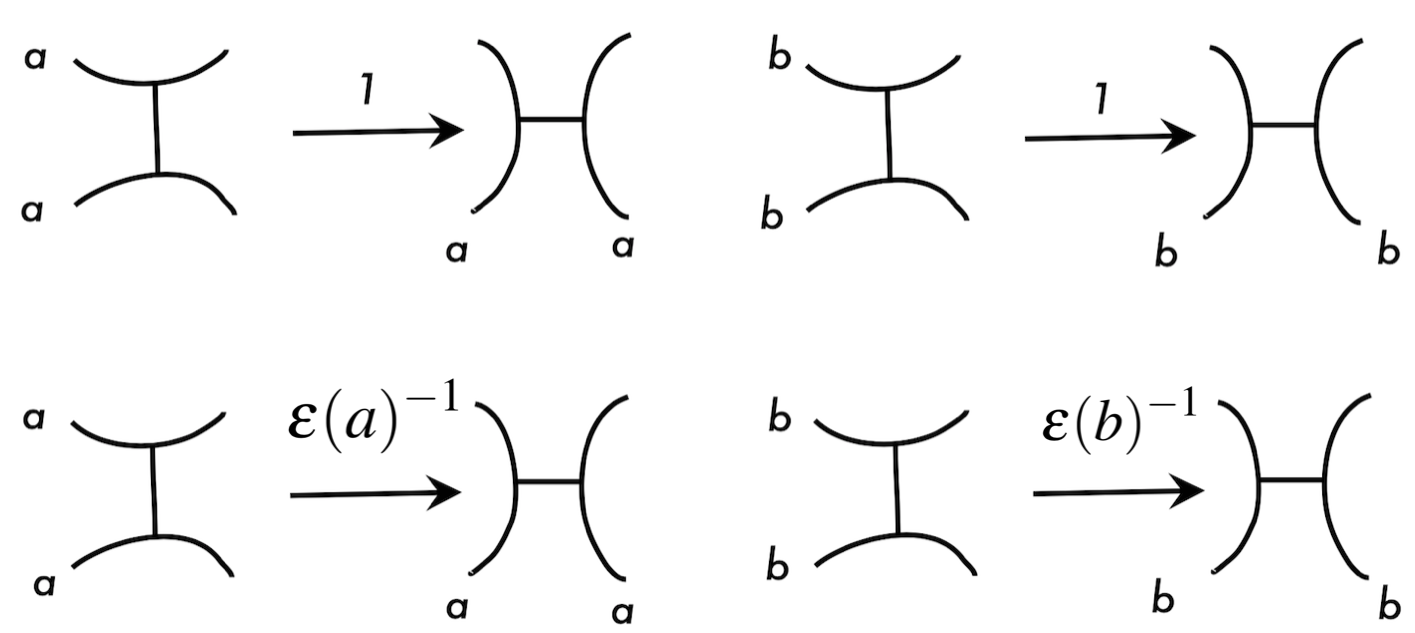}
\end{center}

In all cases the arc-sets and labellings correspond.
 
\end{proof}

Suppose now we are given an enhanced state $(s,x)$  with arc-set $S$. Let  $(s^\prime,x^\prime)$  be an enhanced state in which $s^\prime$ is obtained from $s$ by surgery along arcs in $S$ and $x^\prime$ is the labelling induced by $x$ (using the natural bijection between components of the arc resolution as above). Then
the arc-set of $(s^\prime,x^\prime)$ corresponds to that of  $(s,x)$ and  $(s^\prime,x^\prime)$ is related to  $(s,x)$  in $ \bE$ (the possible cases arising are those presented at the end of the proof of the previous proposition).  

This shows that there is a Boolean sub-poset $\bB(S)$ of $\bE$ containing the enhanced state $(s,x)$ which consists of all enhanced states $(s^\prime,x^\prime)$  obtained from $(s,x)$ by surgery along arcs in $S$ with induced labelling. Proposition \ref{prop:enhancedcompare} shows that there are no other elements of $\bE$ related to $(s,x)$, so $\bB(S)$ consists of {\em all} elements related to $(s,x)$. From this we see that the poset of enhanced states is a disjoint union of Boolean lattices.

\vspace{1.5mm}

\begin{proposition}\label{thm:split}
Let $A$ be a diagonalisable Frobenius algebra. The poset of enhanced states of a link diagram based on $A$ decomposes as a disjoint union of Boolean lattices.
\end{proposition}


\begin{remark}\label{rem:rank} A Boolean component of $\bE$ corresponding to an arc-set $S$ has rank $|S|$.
\end{remark}

\vspace{4mm}

\section{Link homology based on a diagonalisable Frobenius algebra}
Suppose we have a Boolean lattice $B$ equipped with a sign assignment to make squares anti-commute (i.e. of the kind used in the construction of Khovanov homology) and a labelling by modules such that each module is isomorphic to $R$ (the ground ring) and each edge map is an isomorphism. We may follow the usual prescription to build a complex $C^*(B)$.  

\vspace{1.5mm}

\begin{lemma}\label{lem:acyclic}
If rank$(B)>0$, then the complex $C^*(B)$ just defined is acyclic.
\end{lemma}

\vspace{1.5mm}

\begin{proof}
This is clearly true in the rank 1 case and the result is shown by induction using the long exact sequence splitting the Boolean lattice into two co-dimension one faces.

\end{proof}

Returning to the poset of enhanced states of a diagram built from diagonalisable Frobenius algebra with diagonal basis $\{a,b\}$, Proposition \ref{thm:split} shows that the Khovanov complex $C^*$ splits as a direct sum of complexes $C^*(B)$:
$$
C^* = \bigoplus_{B} C^*(B)
$$
where the sum is over Booleans lattices $B$. 

 Lemma \ref{lem:coprod} along with the fact that the multiplication is diagonal ensure that in each such Boolean lattice, all edge-maps are isomorphisms. We can thus apply Lemma \ref{lem:acyclic} to see that only the rank zero Boolean lattices will contribute to the homology. Recalling that the rank of such a Boolean is the number of arcs in the arc-set of any representative we get the following.

\vspace{1.5mm}
 
 \begin{theorem}\label{thm:decomp}
 Let $A$ be a diagonalisable Frobenius algebra defined over a commutative ring $R$. Let $H^*(-)$ be a Khovanov-type link homology constructed from $A$. For a diagram $D$ with poset of enhanced states $\bE$ we have
 $$
 H^*(D) \cong \bigoplus_{(s,x)} R
 $$
 where the sum is over all enhanced states $(s,x)\in \bE$ with empty arc-sets.
 \end{theorem}

\vspace{1.5mm}
 
 The remaining task is to identify those enhanced states having empty arc-sets. Begin by checkerboard colouring the diagram $D$ with the unbounded region coloured white. This induces a black-white colouring of each Kauffman state.  We may now assign an orientation to circles in the (black-white coloured) enhanced state by the following prescription, noting that each circle has a well defined inside and outside:
\begin{center} 
\includegraphics[width=0.9\linewidth]{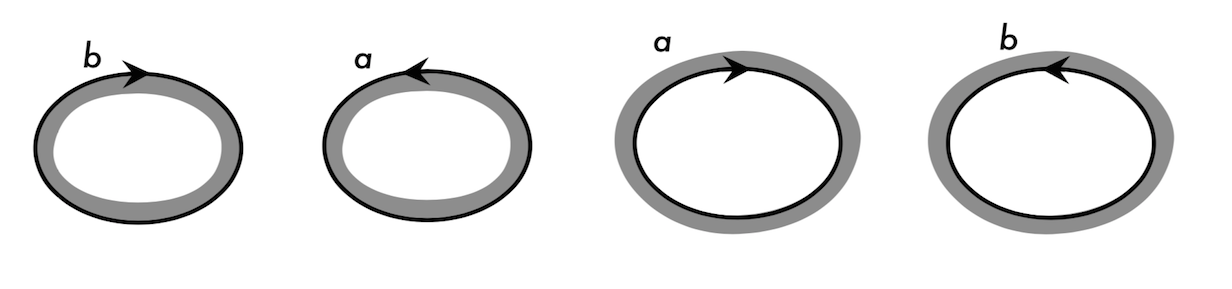}
\end{center}
Remembering the arcs for a moment, near a crossing we have a local piece of one of two types:
\begin{center} 
\includegraphics[width=0.45\linewidth]{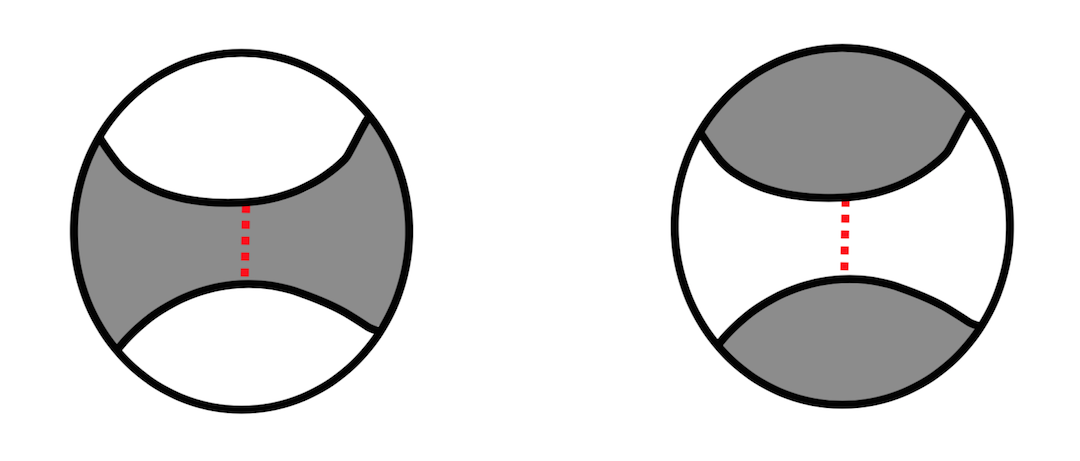}
\end{center}
Since in the arc-resolution under consideration we remove all arcs, this means that in each case the two strands shown are labelled differently. In particular they belong to different circles.

 It easy to verify that all possible ways of closing up the strands gives an orientation for which the two strands in the local piece are oriented the same direction.
It is then possible to re-insert the crossings in a manner that preserves these orientations.
In this way we obtain a map
$$
\{\text{arc-resolutions with no arcs}\} \ra \{\text{orientations of $D$}\}.
$$
This process is evidently invertible: given an orientation of $D$, checkerboard colour as before and then resolve each crossing in the orientation preserving way
This results in a collection of circles and a black-white colouring. Now use the same prescription above to assign $a$ or $b$ to each oriented circle. Indeed this inverse procedure is the method used by Lee to obtain generators from orientations. Summarising:

\vspace{1.5mm}

\begin{proposition}\label{thm:bij}
There is a bijection
$$
\{\text{arc-resolutions with no arcs}\} \longleftrightarrow \{\text{orientations of $D$}\}.
$$
\end{proposition}

\vspace{1.5mm}

Combining Theorem \ref{thm:decomp} and Proposition \ref{thm:bij} yields Theorem \ref{thm:main} as an immediate consequence.

\vspace{4mm}

\subsection*{Acknowledgements} I
 thank Robert Lipshitz and Sucharit Sarkar for providing a sketch argument  and for encouraging me to develop it.

\bibliographystyle{amsalpha}
\bibliography{DiagonalFrobeniusBib}
\end{document}